\documentclass[conference]{IEEEtran}
\usepackage{times}

\usepackage{multicol}
\usepackage[bookmarks=true]{hyperref}
\usepackage[numbers]{natbib}
\newcommand{\tr}{\mathsf{T}}  
\usepackage{xcolor}
\usepackage{amsmath}
\usepackage{amssymb}

\usepackage{graphicx}
\graphicspath{ {./images/} }
\usepackage{caption}    
\usepackage{subcaption} 
\newcommand\norm[1]{\left\lVert#1\right\rVert} 
\usepackage{booktabs} 
\usepackage{makecell} 
\usepackage{siunitx} 

\begin{document}

\title{Optimal Control of Granular Material}

\author{
\authorblockN{Yuichiro Aoyama}
\authorblockA{School of Aerospace Engineering\\
Georgia Institute of Technology\\
Atlanta, Georgia\\
yaoyama3@gatech.edu}
\and
\authorblockN{Amin Haeri}
\authorblockA{University of Chicago\\
Chicago, Illinois\\
haeri@uchicago.edu}
\and
\authorblockN{Evangelos A. 
 Theodorou}
\authorblockA{School of Aerospace Engineering\\
Georgia Institute of Technology\\
Atlanta, Georgia\\
evangelos.theodorou@gatech.edu}}


%

\maketitle
\begin{abstract}
The control of granular materials, showing up in many industrial applications, is a challenging open research problem. Granular material systems are complex-behavior (as they could have solid-, fluid-, and gas-like behaviors) and high-dimensional (as they could have many grains/particles with at least 3 DOF in 3D) systems. Recently, a machine learning-based Graph Neural Network (GNN) simulator has been proposed to learn the underlying dynamics. In this paper, we perform an optimal control of a rigid body-driven granular material system whose dynamics is learned by a GNN model trained by reduced data generated via a physics-based simulator and Principal Component Analysis (PCA).  We use Differential Dynamic Programming (DDP) to obtain the optimal control commands that can form granular particles into a target shape. The model and results are shown to be relatively fast and accurate. The control commands are also applied to the ground-truth model,(i.e., physics-based simulator) to further validate the approach.
\end{abstract}
\IEEEpeerreviewmaketitle

\section{Introduction}
The control of granular materials has many applications such as mining machinery, industrial robots (e.g. excavators), and robots that work with humans \cite{Hemami1994loader, Connor2017learninggranular, Mateo2019GaussainGranular}. There are various aspects to consider for solving such a complex problem. One difficulty here is obtaining the granular material model. Another difficulty is that granular materials are high-dimensional systems, as they could include many grains/particles with at least 3 DOF (in 3D). To achieve the system model, machine learning methods have been recently shown as efficient alternatives to traditional numerical methods. They require a huge amount of training data that should be generated either by verified computer simulations or experiments (which are usually not practically feasible). Moreover, the control of such systems by actively interacting with the particles has not been addressed in literature yet. We will discuss it in more detail, in the following. 

\textbf{Simulation Method.}
Physics-based numerical granular flow simulation methods can be considered highly accurate candidates as the system model \cite{Haeri20Gravity}. But they predominantly suffer from their intensive computational complexity, especially in hardware-constrained robotics applications \cite{Haeri20isarc}. However, Machine learning (ML) has begun to be applied to the problem of reducing the runtime of simulating complex physics \cite{Haeri21Aero}. Data-driven ML models are an emerging class of physics simulation methods. They are trained by either computer-generated synthetic or real experimental data. They could be accurate enough depending on the utilized data and model architecture.

\textbf{Learning Model.}
Here, the goal is to properly learn the underlying dynamics principles in the physical systems automatically. Several recent methods have been utilized for fluid simulations, including regression forest, multi-layer perceptron (MLP), autoencoder (AE) convolutional neural network (CNN), generative adversarial network (GAN), and loss function-based method \citep{Lad17regFor, Wie19ff, Thu20ae, Kim19cnn, Xie18gan, Pra20lf}. They can mostly handle a specific state of materials (either gas-, solid-, or specifically fluid-like behaviors). However, graph neural networks (GNN) \citep{FrancoGNN2009} are generalized models proven to be capable of simulating different states of materials. Recent state-of-the-art research has developed a GNN model \citep{sanchez2020learning}, which has outperformed some other recent methods \citep{li2019learning, Umm20cc} in terms of accuracy and ease of implementation. However, aside from the time complexity, GNNs might be infeasible for desktop GPUs with mid-level limits on memory (i.e. VRAM of 8-20GB). 

\begin{figure}[!b]
    \centering
    \includegraphics[width=0.98\linewidth]{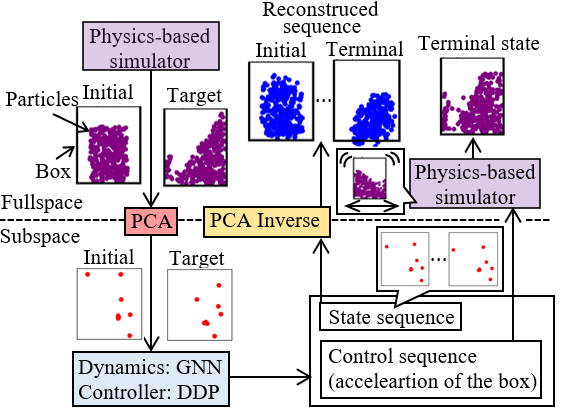}
    \caption{Our model overview.}\label{fig:schematic}
\end{figure}

\textbf{Reduced Model.}
In physics simulations, model reduction methods are utilized to capture the effective degrees of freedom of a physical system \citep{subs1,subs2,subs3,subs4}. This is to reduce both processing time and memory space usage corresponding to the system. Recent literature \citep{haeri2021subspace, pcamg} has shown that principal component analysis (PCA) \cite{pearson1901PCA} could be an effective choice to capture the primary modes of object deformations while keeping the expressive power of the ML model. PCA can significantly reduce the training and rollout runtime and memory usage of ML models by learning the reduced dynamics (i.e. subspace learning \citep{SINDy}).

In \cite{haeri2021subspace}, in order to run the ML-based simulation in real-time, a subspace GNN model has been trained by a reduced number of particles obtained via PCA. The subspace output of the GNN model can then be projected back to fullspace by the PCA inverse. This approach has been able to accurately predict the interactions of rigid bodies and granular flows (e.g. the interactions of an excavator bucket or a rover wheel with soil particles). Note that the training data has been generated by an accurate and efficient continuum physics-based simulation method called Material Point Method (MPM) with a nonlocal granular flow constitutive model (NGF) \cite{haeri2022three, Kamrin2019nonlocal}. This method can also predict the interaction forces applied to the rigid body by the granular material.

\textbf{Particle Control.} 
There are a few granular material dynamics models proposed for control purposes in literature. In \cite{Tuomainen2022particle}, using a particle dynamics model introduced in \cite{sanchez2020learning}, a control problem, where particles are shaped to be a target shape, is solved. The task is to find a suitable place to pour particles (instead of actively interacting with them as done in this paper). Also, in \cite{li2019learning}, some control problem examples including forming fluids and deformable material into target shapes, are presented. However, the particles are assumed to not have real mechanical properties which limits their potential real-world applications. 

Control of particle systems can be seen as a trajectory optimization problem where state and control sequences are optimized to reach a desired state (i.e. a target shape). For trajectory optimization, Sequential Quadratic Programming (SQP)\cite{Boggs1995SQP} and Differential Dynamic Programming (DDP)\cite{Jacobson1970ddp} have been widely used in robotic applications. The system state can be represented by a concatenation of particle states and actuators (e.g. bucket and robotic arm). Thus, the control state dimension is significantly smaller than that of the system state dimension. For this type of system, DDP is computationally noticeably faster than SQP. This is because in DDP, size of Hessian is determined by the control state dimension (whereas the size of Hessian in SQP is determined by the sum of the system and control state dimensions) \cite{Liao1991convDDP, GILL2000SQP}.

\textbf{Contribution.}
In this paper, we aim to actively control the shape of a granular material interacting with a rigid body and to form its shape into a target shape. We perform an optimal control using the DDP algorithm with the PCA-reduced dynamics learned by the GNN model \cite{haeri2021subspace}. We use a 2D rectangular rigid box with soil inside as an example system whose control inputs are the box temporal accelerations. First, we generate a training dataset using a physics-based simulator \cite{haeri2022three} including several examples each with thousands of particles and different rigid box properties. Second, the granular flow system is reduced via PCA (i.e. from ~1000 particles to 8 modes). We then train a GNN model, obtaining system dynamics in subspace. By applying DDP to the reduced dynamics, the box-related control input commands are obtained. The control commands can cause the  particles to form a target shape. Finally, we apply the control commands to the physics-based simulator to further validate our approach. An overview of our model is shown in Fig. \ref{fig:schematic}. To our knowledge, such a framework has not already been proposed for the granular flow control problem.

\section{Preliminaries}
In this section, we will provide a background on the major components used in our framework including the Material Point Method (MPM) physics-based simulation method, the machine learning-based Graph Neural Network (GNN) simulation model, the Principal Component Analysis (PCA) model reduction method, and the Differential Dynamic Programming control algorithm.

\subsection{Physics-Based Numerical Simulator}
\label{subsec:Physics-based}
Material Point Method (MPM) with a Nonlocal Granular Fluidity constitutive model (NGF) is a continuum numerical simulation method specifically developed for granular material simulations. It is used to solve the underlying partial differential equations. This method can simulate the interactions of rigid bodies with granular flows accurately and efficiently (not real-time). In MPM, rigid bodies and granular materials are represented by particles of different types. Also, the motion of rigid bodies could be predefined. In other words, one can input the rigid body trajectories into MPM. Here, this method is used to generate training data. We refer the reader to \cite{haeri2022three} for detailed derivations and procedures.

\subsection{Graph Neural Network Simulator}
\label{sec:GNN}
The Graph Neural Network (GNN) simulator utilizes a graph representation of particles possibly with different types (e.g. rigid, granular, etc.) \cite{sanchez2020learning}. Let us define a graph $G = (V, E)$, where $v_{i} \in V$ and $e_{i,j} \in E$ are the node and edge attributes, respectively. The subscripts $i$, $j$ in $e_{i,j}$ indicate the sender and receiver node indices. In the graph representation of particles, nodes $v_{i}$ correspond to particles and edges $e_{i,j}$ correspond to the relations between the particles. More precisely, if two particles are closer than a threshold, they can be connected with edges with specific attributes. Note that we do not aim to use such relations as in subspace the relations might not be proximity-based.

\begin{figure}[!t]
    \centering
    \includegraphics[width=0.98\linewidth]{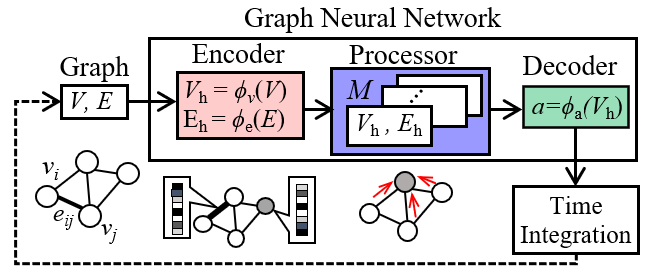}
    \caption{Graph Neural Network (GNN) simulator.}\label{fig:GNS}
\end{figure}

Fig. \ref{fig:GNS} shows the overall procedure performed in the GNN simulator. It takes in a graph and outputs the particle accelerations. The node and edge attributes are given by
\begin{align}
    v_{i} &= [q^{i}_{k},\underbrace{\dot{q}^{i}_{k-C+1},\cdots, \dot{q}^{i}_{k}}_{C \text{ sequence} }, s^{i}, g_{k}], \\
    e_{i,j} &= [q^{i}_{k}-q^{j}_{k}, \norm{q^{i}_{k}-q^{j}_{k}}], \nonumber
\end{align}
where $q^{i}_{k}$ is the position of $i$ th particle at time step $k$, $s^{i}$ is particle type (e.g., sand, rigid body), and $g_{k}$ is global feature (shared among all nodes). Note that $C$ sequence of velocity can be obtained by $C+1$ sequence of positions using finite difference schemes. GNN first encodes an input graph into a latent graph via encoder Multi-Layer Perceptrons (MLPs) denoted by $\phi_{v}$ and $\phi_{e}$ as
\begin{align}
    V_{h} = \phi_{v}(V), \quad E_{h} = \phi_{e}(E).
\end{align}
Then, the processor performs $M$ times of message-passing steps to propagate information and update the node and edge attributes. The detailed equations of the processor component can be found in \citep{haeri2021subspace}.
Finally, the particle accelerations are predicted via a decoder MLP $\phi_{\rm{a}}$ as follows
\begin{align}
    \ddot{q}_{k+1}^{i} = \phi_{\rm{a}}(v_{i}).
\end{align}
The accelerations are used to update node positions $q_{k}$ using implicit Euler integration with time step $\Delta t$ via
\begin{align}
    \dot{q}_{k+1} = \dot{q}_{k} + \ddot{q}_{k}\Delta t, \ {q}_{k+1} = {q}_{k} + \dot{q}_{k+1}\Delta t  
\end{align}
where we drop $i$. $q_{k+1}$ is used as the node positions of the new input graph to forward the simulation. Note that similar to the physics-based simulator, the motion of rigid particles can be predefined. Therefore, the GNN model can be learned to predict the motion of normal granular particles only.
 
\subsection{Principle Component Analysis} \label{subsec:pre_PCA}
Principal Component Analysis (PCA) is a non-parametric linear dimensionality reduction method \cite{Hotelling1933AnalysisOA}. It provides a data-driven, hierarchical coordinate system to re-express high-dimensional correlated data. The resulting coordinate system geometry is determined by principal components (also known as modes). These modes are uncorrelated (orthogonal) to each other but have a maximal correlation with the observations. Particularly, the goal of this method is to compute a subset of ranked principal components to summarize the high-dimensional data while retaining trends and patterns.

Here, we explain how PCA is used to obtain the reduced data. Define $n_{\rm{n}}$ as the number of normal particles in each example and $m_{\rm{o}}$ as the number of observations. Note that we are not using rigid particles. Given that the system is $d$ dimensional with time horizon $N$, and the number of training examples is $n_{\rm{e}}$, the number of observations is given by $m_{\rm{o}} = n_{\rm{e}} \times N$. In 2D with $d=2$, the data matrix $X_{\rm{d=2}}$ of the particle trajectories (i.e. positions) is written as
\begin{align}\label{eq:PCA_data_matrix}
    X_{\rm{d}} = \begin{bmatrix}
    x_{1}^{1} & x_{1}^{2}&  \cdots & x_{1}^{n_{\rm{n}}}\\
    y_{1}^{1} & y_{2}^{2}&  \cdots & y_{1}^{n_{\rm{n}}}\\
              &           & \vdots &                \\
    x_{m_{\rm{o}}}^{1} & x_{m_{\rm{o}}}^{2}& \cdots & x_{m_{\rm{o}}}^{n_{\rm{n}}}\\
    y_{m_{\rm{o}}}^{1} & y_{m_{\rm{o}}}^{2}& \cdots & y_{m_{\rm{o}}}^{n_{\rm{n}}}\\ 
    \end{bmatrix}
    \in \mathbb{R}^{dm_{\rm{o}} \times n_{\rm{n}}},
\end{align}
where superscripts and subscripts of $x$ and $y$ are indices for the particle and observation number, respectively. Notice that the column of $X_{\rm{d}}$ represents a flattened and concatenated trajectory of a particle, whereas the row represents the position of all the particles at the same time step.
In order to perform PCA, we first compute the mean $\mu$ and covariance matrix $S$ of the observations via
\begin{align*}
    \mu &= \frac{1}{m_{\rm{o}}}\sum_{i=1}^{m_{o}}[X_{\rm{d}}]_{i,:}\in \mathbb{R}^{1 \times n_{\rm{n}}} \\
    S &= \frac{1}{m_{\rm{o}}}(X_{\rm{d}}-\mathbf{1}\mu)^{\mathsf{T}}(X_{\rm{d}}-\mathbf{1}\mu) \in \mathbb{R}^{n_{\rm{n}}\times n_{\rm{n}}}
\end{align*}
where $[X_{\rm{d}}]_{i,:}$ is the $i$'th row of $S$ and $\mathbf{1} = [1, \cdots, 1]^{\mathsf{T}} \in \mathbb{R}^{n_{p}}$.
$n_{\rm{nr}}$ principle components of centralized (i.e. zero-mean) data are given as eigenvectors of the covariance matrix of data, which corresponds to the first $n_{\rm{nr}}$'th largest eigenvalues, forming the loading matrix $W$ as
\begin{align}
    \label{eq:PCA W}
    W = [w_{1}, \cdots w_{n_{\rm{_pr}}}] \in \mathbb{R}^{n_{\rm{n}} \times n_{\rm{nr}}}.
\end{align}
Multiplying $W$ from the right, the centralized data is mapped to the subspace as
\begin{align}\label{eq:PCA projection}
    Z = [X_{\rm{d}}-\mathbf{1}\mu]W \in \mathbb{R}^{dm_{\rm{o}}\times{n_{\rm{nr}}}}.
\end{align}
where $\mu$ is the mean of data. Notice that the number of particles is now reduced from $n_{\rm{n}}$ to $n_{\rm{nr}}$. $Z$ can be mapped back to fullspace, resulting in the reconstruction of $X_{\rm{d}}$ denoted by $X^{\dagger}_{\rm{d}}$ 
\begin{align}
    \label{eq:PCA reconstruction}
    X_{\rm{d}}^{\dagger} = ZW^{\mathsf{T}} + \mathbf{1}\mu.
\end{align}
When $n_{\rm{n}} = n_{\rm{nr}}$ and $W^{\mathsf{T}} = W^{-1}$, the $X_{d}$ is perfectly reconstructed.

\subsection{Differential Dynamics Programming}
Consider the discrete-time optimal control problem of a dynamical system given by
\begin{align}\label{eq:unconstrained_optimalcontrol}
    \min_{{U}}&\hspace{0.8mm}J({X},{U})=\min_{{U}}\big[{\sum_{k=0}^{N-1} l({x}_k,{u}_k)}+\Phi({x}_N)\big] \\
    &\text{subject to}\hspace{3mm} {x}_{k+1}={f}({x}_k,{u}_k),\quad k=1,2,...,N-1 \nonumber
\end{align}
where ${x}_k\in\mathbb{R}^n$, ${u}_k\in\mathbb{R}^m$ denote the state and control inputs of the system at time step $t_k$, respectively. ${f}:\mathbb{R}^n\times \mathbb{R}^m\rightarrow\mathbb{R}^n$ is the transition dynamics function. The scalar-valued functions $l(\cdot,\cdot)$, $\Phi(\cdot)$, and $J(\cdot)$ denote the running, terminal, and total cost of the problem, respectively. We also let ${X} =[{x}_{0},\dots,{x}_{N}]$, ${U} =[{u}_{0},\dots,{u}_{N-1}]$ be the state and control trajectory over the time horizon $N$. The cost-to-go at time step $k=i$, i.e. cost starting from $k=i$ to $N$ is given by
\begin{align}\label{eq:DDP_cost_to_go}
    {J}_{i}(X_{i},U_{i}) = \big[\sum_{k=i}^{N-1} l(x_k,u_k)\big] + \Phi(x_{N})
\end{align}
with $X_{i} =[{x}_{i},\dots,{x}_{N}]$, $U_{i} =[{u}_{i},\dots,{u}_{N-1}]$. The value function is defined as the minimum cost-to-go at each state and time step via
\begin{equation}
    \label{value-function-definition}
    V_k({x}_k):=\min_{{u}_k}J({X},{U}).
\end{equation}

Given Bellman's principle of optimality that provides the following rule
\begin{equation}\label{bellman}
    V_k(x_k) =  \min_{{u}_k} [l({x}_k,{u}_k) + V_{k+1}({x}_{k+1})],
\end{equation}
DDP finds the local solutions to Eq. \eqref{eq:unconstrained_optimalcontrol} by expanding both sides of the rule (i.e. Eq. \eqref{bellman}) about given nominal trajectories, $\bar{X}$ and $\bar{U}$. To examine this, we define $Q$ function as the argument of $\min$ on the right-hand side of Eq. \eqref{bellman}:  
\begin{equation}
    \label{Qfunction}
    Q_k(x_k,u_k) = l(x_k,u_k) + V_{k+1}(x_{k+1}).
\end{equation}
Taking quadratic expansions of $Q_k$ about $\bar{X}$ and $\bar{U}$ by considering deviation $\delta x_k=x_k -\bar{x}_{k}$, $\delta u_k=u_{k} - \bar{u}_{k}$,
and mapping terms in both sides of Eq. \eqref{Qfunction} with linearized dynamics $\delta x_{k+1} =  f_{x} \delta x_{k} + f_{u} \delta u_{k}$ result in the $Q$ derivatives (evaluated on $\bar{X}$ and $\bar{U}$) as follows
\begin{equation}
    \label{Qexpanded}
    \begin{split}
        &{Q}_{{xx},k} ={l}_{{xx}}+{f_x}^{\tr}{V}_{{xx},k+1}{f_x},\hspace{0.8mm}{Q}_{{x},k}={l}_{{x}}+{f_x}^{\tr}{V}_{{x},k+1}
         \\
        &{Q}_{{uu},k} = {l_{{uu}}}+{f_u}^{\tr}{{V}_{{xx},k+1}}{f_u}, \hspace{0.8mm}{Q}_{{u},k}={l_{ u}}+{f_u}^{\tr}{{V}_{{x},k+1}}
        \\
        &{Q}_{{xu},k} ={l_{{xu}}}+{f_x}^{\tr}{{V}_{{xx},k+1}}{f_u}.
    \end{split}
\end{equation}

By considering quadratically expanded $Q$, we can explicitly optimize with respect to $\delta{u}_{k}$ and compute the locally optimal control update via
\begin{align}\label{eq:delta-u-star}
    & \delta {u}^{\ast}_k={\kappa}+{K}\delta {x}_{k} \\
    \text{with}\quad &\kappa:= -{Q}^{-1}_{{uu}}{Q_{{u}}},\hspace{1.8mm} {K} = -{Q}^{-1}_{{uu}}{Q_{{ux}}} \nonumber
\end{align}
where $\kappa$ and $K$ are known as feedforward and feedback gains, respectively.
Note that we drop time indices for $Q$ for readability. Now, $\delta {u}_{k}^{\ast}$ is computed using $V_{x,k+1}$ and $V_{xx,k+1}$. Notice the difference in time instances of $u_{k}$ and $V_{x,k+1}$. To propagate $V_{k}$ back in time, we plug quadratic expansion of $V_{k+1}$ with linearized dynamics and $\delta {u}^\ast_k$ into Eq. \eqref{bellman}, and comparing coefficients of $\delta x_{k}$, obtaining
\begin{align} \label{eq:value_riccati}
    V_{{x},k}&={Q_{{x},k}}-{Q_{{xu},k}}{Q_{{uu},k}^{-1}}{Q_{{u},k}},\\
    V_{{xx},k}&={Q_{{xx},k}}-{Q_{{xu},k}}Q_{{uu},k}^{-1}{Q_{{ux},k}}, \nonumber
\end{align}
with boundary condition $V_{N} = \Phi(x_{N})$. The process of computing $V_{x}$ and $V_{xx}$ is called backward pass. After the completion of the backward pass, a new state-control sequence is determined by propagating the dynamics forward in time (called forward pass).

Typically, a backtracking line-search with step size $\alpha \in (0,1]$ is performed as follows
\begin{align}
    \bar{u}_{k}^{\rm{new}} &= \bar{u}_{k} + \alpha k + K\delta x_{k}, \quad \delta x_{k} =  x_{k}^{\rm{new}}- x_{k} \\\bar{x}_{k+1}^{\rm{new}} &= f(\bar{x}^{\rm{new}}_{k},\bar{u}^{\rm{new}}_{k}). \nonumber
\end{align}
After finding $\alpha$ that reduces the cost, the new sequence is treated as the new nominal trajectory for the next iteration. The procedure of forward and backward passes is repeated until certain convergence criteria are satisfied.

To ensure convergence, $Q_{{uu}}$ must be regularized when its positive definiteness cannot be guaranteed \cite{Liao1991convDDP} i.e.
\begin{align*}
    Q_{uu}^{\rm{reg}} = Q_{uu} + \lambda I_{m}
\end{align*}
with $\lambda>0$ and an identity matrix $I$. Practical scheduling of $\lambda$ is found in \cite{TassaDDP2012}.
To compute the gains, $Q_{uu}^{\rm{reg}}$ is replaced with $Q_{uu}$ in Eq. \eqref{eq:delta-u-star}. With regularized $Q_{uu}$ and gains, Eq. \eqref{eq:value_riccati} can be written as
\begin{align*}
    V_{x,k} &= Q_{{x}} + {K}^{T}Q_{{uu}}{\kappa}
    +{K}^{T}Q_{{u}} + Q_{{ux}}^{T}{\kappa} \\\notag
    V_{xx,k}&= Q_{{xx}} + {K}^{T}Q_{{uu}}{K}
    +{K}^{T}Q_{{ux}} + Q_{{ux}}^{T}{K}
\end{align*}
where regularized $Q$ is used in $\kappa$ and $K$. Finally, we note that this derivation uses the first-order expansion of dynamics, whereas the original DDP uses the second-order expansion. We use the first-order expansion due to its computational efficiency and numerical stability \cite{iLQR2004Li}.

\section{Learned Dynamics and Optimal Control}
The dynamics of our system, represented by a GNN model, is explained here. Our full system is defined as a rectangular box containing granular material. Furthermore, the system state-space representation is derived to feed into the DDP method.

\subsection{System Dynamics}
Let $p_{k}$ be the positions of all the particles including normal particles and rigid particles which represents the box at time step $k$. Note that we distinguish normal particles and rigid particles later. Also, only normal particles are projected on to subspace via PCA. Let $b_{k}$ be a representative position of the box, which pick as the center of the bottom side of the box. The particle positions $p_{k}$ are given in the box frame, whereas the box positions $b_{k}$ are in the world frame.

Concatenating $p_{k}$ and $b_{k}$, we define 
\begin{align}
    x_{k} = [p_{k}^{{\tr}}, b_{k}^{{\tr}}]^{\tr}.
\end{align} 
We define vector $X$ as a sequence of the concatenated particle positions via 
\begin{align}
    \label{eq:DDP X state}
    X_{k:k+C} = [p_{k}^{\tr},b_{k}^{\tr},\cdots p_{k+C}^{\tr},b_{k+C}^{\tr}]^{\tr}\in \mathbb{R}^{d(n_{p}+1)(C+1)}
\end{align}
where $C$ is the length of the sequence of velocity introduced in \S \ref{sec:GNN}. The system's control input $u_{k}$ is the acceleration of the box. $C+1$ sequence of position is required to recover $C$ sequence of velocity. The GNN model takes in $X_{k:k+C}$ as inputs, and predicts normal particle accelerations $\ddot{p}_{k+C+1}$ as detailed in \S \ref{sec:GNN}. In addition to $X_{k:k+C}$, we add the position, velocity, and acceleration of the box as global features which are shared among all the particles, and thus nodes of the GNN. $u_{k+C}$ here updates rigid particles.

Let us write only particle position part of $X_{k:k+C}$ as $P_{k:k+C}$, i.e., $P_{k:k+C} = [p_{k}^{\tr},\cdots p_{k+C}^{\tr}]^{\tr}$. Similarly, let us write the box position part as $B_{k:k+C} = [b_{k}^{\tr},\cdots b_{k+C}^{\tr}]^{\tr}$. By using the finite difference method and Euler integration, the relationship is compactly written as  
\begin{align}
    p_{k+C+1} &= f_{\rm{p}}(P_{k+C}, b_{k+C}, \dot{b}_{k+C},\ddot{b}_{k+C}, u_{k+C})\\&= f_{\rm{p}}(P_{k:k+C}, B_{k:k+C}, u_{k+C}). \nonumber
\end{align}
In the following, it will be explained in detail.

\subsection{Normal and Rigid Particles}
Let $g_{\rm{GNN}}(\cdot)$ be a function of the trained GNN simulator. We obtain the accelerations of the most recent element of the state sequence by feeding the sequence of positions
\begin{align}
    \ddot{p}_{k+C+1} = g_{\rm{GNN}}(X_{k:k+C}).
\end{align}
Using this acceleration, update velocity and position using the semi-implicit Euler scheme;
\begin{align}
    \dot{p}_{k+C+1} &= {\dot{p}_{k+C}} + \ddot{p}_{k+C+1}\Delta t\\ 
    {p_{k+C+1}} &= p_{k+C} + \dot{p}_{k+C+1} \Delta t  \nonumber
\end{align}
As explained in \S \ref{sec:GNN}, only normal (not rigid) particles are updated by GNN. Since the particles are represented in the box coordinates, the rigid particles do not change their positions in the box frame. Therefore,
\begin{align}
    [p_{k+C+1}]^{\rm{b}} &= [p_{k+C}]^{\rm{b}}
\end{align}
where $[p_{k+C}]^{\rm{b}}$ stands for rigid particles for the box.

\subsection{Box Position}
As mentioned, the box acceleration is the control input, $\tilde{u}_{k+C}$. Hence, the box velocity and position are updated via
\begin{align}
    \dot{b}_{k+C+1} &= {\dot{b}_{k+C}} + \tilde{u}_{k+C} \Delta t\\
    b_{k+C+1} &= b_{k+C} + \dot{b}_{k+C+1} \Delta t \nonumber
\end{align}
where $\tilde{u}_{k+C}$ distributes $u_{k+C}$ to the corresponding dimensions. In our 2D example, $b_{k} \in \mathbb{R}^{2}$, whereas $u_{k} \in \mathbb{R}$ and thus they cannot be added as they are.
\begin{align}
    \tilde{u}_{k} = [u_{k},0]^{\tr}.
\end{align}

\subsection{State-Space Representation}
Given all the aforementioned updates, the dynamics of the system is written in a matrix-vector form via
\begin{align}\label{eq:state_space_matrix}
    &X_{k+C+1} =
    \begin{bmatrix}
        p_{k+1} \\b_{k+1} \\ \vdots \\p_{k+C} \\ b_{k+C} \\ p^{\rm{n}}_{k+C+1} \\
        p^{\rm{b}}_{k+C+1} \\
        b_{k+C+1} 
    \end{bmatrix}
    =
    \begin{bmatrix}
        M
        \begin{bmatrix}
            p_{k} \\ b_{k}\\ \vdots \\ p_{k+C} \\  b_{k+C} 
        \end{bmatrix}\\
        [g_{\rm{GNN}}(P_{k+C},B_{k+C})]^{\rm{n}} \\
        p_{k+C}^{\rm{b}}\\
        2b_{k+C} - b_{k+C-1} + \tilde{u}_{k+C}(\Delta t)^{2}
    \end{bmatrix}\\\notag
&{\text{where}}\\\notag
    &M = [O_{d(n_{\rm{p}}+1)C, d(n_{p}+1)}, M_{0}]\in \mathbb{R}^{d(n_{{\rm{p}}+1})C \times d(n_{\rm{p}}+1)(C+1)}\\\notag
    &M_{0} = {\rm{blkdiag}}[I_{dn_{\rm{p}}}, I_{d}, \cdots, I_{dn_{\rm{p}}}, I_{d}]\in\mathbb{R}^{d(n_{\rm{p}}+1)C\times d(n_{\rm{p}}+1)C}
\end{align}
and $[g_{\rm{GNN}}]^{\rm{n}}$ is the output of the GNN model for normal particles. $O$ and $I$ are zero and identity matrices, respectively, whose sizes are specified in their subscripts.

From Eq. \eqref{eq:state_space_matrix}, we obtain state-space representation of the system dynamics as
\begin{align}\label{eq:box_particle_dyn_state_space_rep}
    X_{k+C+1} = F(X_{k+C}, u_{k+C}).
\end{align}

\subsection{Jacobian Matrix}
Since most parts of the dynamics are linear mappings, the Jacobian matrix of the dynamics is given by 
\begin{align}
    \frac{\partial F}{\partial X_{k'}} 
    &= \begin{bmatrix}
    & M &\\
    & \frac{\partial{[g_{\rm{GNN}}]^{\rm{n}}}}{\partial X_{k'}} & \\
    O_{1} & O_{2} & O_{3} & I_{dn_{\rm{b}}} & O_{2}\\
    O_{4} & - I_{d} & O_{5}& O_{6} & 2I_{d}
    \end{bmatrix},\label{eq:jacobian_matrix_x}\\
    \frac{\partial{F}}{\partial{u_{k'}}} &= \begin{bmatrix}
    O\\
    \frac{\partial \tilde{u}_{k'}}{\partial u_{k'}}
    \end{bmatrix},\label{eq:jacobian_matrix_u}\\\notag
    {\text{with}} \
    O_{1} &= O_{dn_{\rm{b}},d(n_{\rm{p}}+1)(C-1)+n_{\rm{p}}}, \quad O_{2} = O_{dn_{\rm{b}},d} \\\notag
    O_{3} &= O_{dn_{\rm{b}}, dn_{\rm{n}}}, \quad  O_{4} = O_{d,d(n_{\rm{p}}+1)(C-1)+n_{\rm{p}}} \\\notag
    O_{5} &= O_{d,dn_{\rm{n}}}, \quad O_{6} = O_{d,dn_{\rm{b}}},
\end{align}
and $n_{\rm{b}}$ is the number of rigid particles which represents the box.
By utilizing this structure, the Jacobian matrix can be efficiently computed.

\subsection{Optimal Control}
Now, we can set up an optimal control problem (i.e. Eq. \eqref{eq:unconstrained_optimalcontrol}) for our system. The dynamics is given by Eq. \eqref{eq:box_particle_dyn_state_space_rep} with the Jacobian matrix derived in Eq. \eqref{eq:jacobian_matrix_x} and \eqref{eq:jacobian_matrix_u}. We set quadratic running and terminal costs as follows
\begin{align}\label{eq:numerical_example_quad_cost}
    l(X_{k},u_{k}) &= 0.5 u_{k}^{\tr}Ru_{k} + 0.5 (X_{k}-X_{\rm{g}})^{\tr}  P(x_{k}-x_{\rm{g}}), \\ \Phi(X_{N}) &= 0.5(X_{N}-X_{\rm{g}})^{T}Q(X_{N}-X_{\rm{g}}),
\end{align}
with desired state $X_{g}$ and positive semi-definite weight matrices $P,Q$, and $R$. 
Since we now have everything required for DDP, we can compute the optimal control sequence for our task, i.e., the sequence of box accelerations to form the particles into a target shape.

\section{Numerical Experiments}
In this section, we will provide the results of our numerical experiments and explain the procedures used to obtain the results. The procedures include data generation, dimensionality reduction, model training, and setting initial and target conditions for DDP.

\subsection{MPM Data}
We generate the training data by shaking a box with granular particles inside using the physics-based simulation method, MPM-NGF, introduced in \S\ref{subsec:Physics-based}. In the dataset, we have generated 324 examples with varying widths, motion frequencies, and motion amplitudes of the box. The granular material properties remain constant through the examples. The data acquisition rate is set to 60Hz, which is also the frequency of the optimal control problem we solve by DDP later on. Also, each example lasts 1.5 seconds (i.e. 90 frames).

In examples with low-frequency box motions ($\lesssim 0.5$ Hz), particles only move with the box and do not show interesting behavior for training purposes. Thus, we aim to set high box motion frequencies to obtain meaningful data for training. However, there exists a stability (also known as Courant) condition for the convergence of the partial differential equations solved in the physics-based simulation method, MPM-NGF \cite{Courant1967}. The condition requires $v \Delta t / \Delta x \leq 1$ with grid size $\Delta x$, time step $\Delta t$, and particle velocity magnitude $v$. With fixed $\Delta x$ and large $v$, it requires a small $\Delta t$, leading to a large computational runtime. We set the maximum frequency of the box to 3Hz by considering a balance between the computational runtime (for the data generation) and the motion of particles empirically.

\subsection{PCA Reduction}
\begin{figure}[!t]
     \centering
     \begin{subfigure}[b]{\linewidth}
         \centering
         \includegraphics[width=\textwidth]{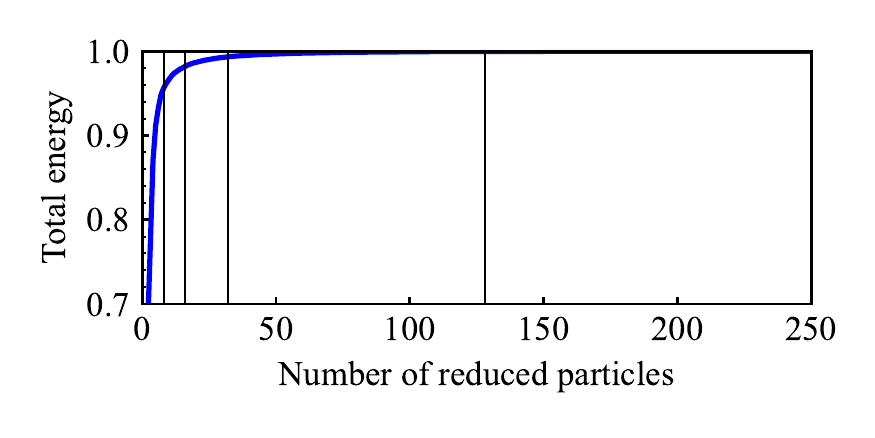}
         \caption{Total energy (cumulative) versus the number of reduced particles (i.e. PCA modes). The vertical lines correspond to 8, 16, 32, 64, and 128 numbers.}
         \label{fig:PCA energy}
     \end{subfigure}
     \hfill
     \begin{subfigure}[b]{0.51\linewidth}
         \centering
         \includegraphics[trim={0.5cm 0.3cm 0.3cm 0.3cm},clip,width=\textwidth]{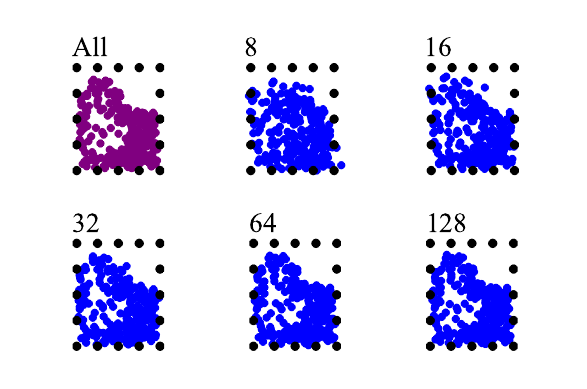}
         \caption{Example 1}
         \label{fig:reconst_ex_narrow}
     \end{subfigure}
     \begin{subfigure}[b]{0.42\linewidth}
         \centering
         \includegraphics[trim={0.5cm 0.3cm 0.8cm 0.3cm},width=\textwidth]{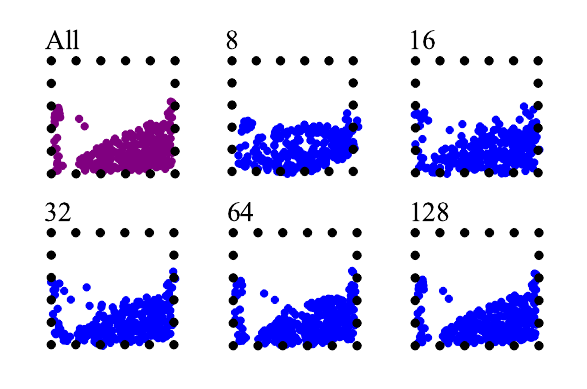}
         \caption{Example 2}
         \label{fig:reconst_ex_wide}
     \end{subfigure}
    \caption{Total energy (a) and reconstruction from reduced particles for two examples (b and c). Particles in fullspace are drawn in purple and reconstructed particles in blue.}
    \label{fig:three graphs}
\end{figure}

After generating the training data, PCA is applied to reduce the dimensionality of the data as explained in \S\ref{subsec:pre_PCA}. Fig. \ref{fig:PCA energy}
shows the relationship between the number of reduced particles (i.e. PCA modes) and the cumulative total energy (i.e. variance). The total energy of the first $i$ reduced particles is given by $\textstyle{\sum_{k=1}^{i}{\lambda_{k}}}/\textstyle{\sum_{k=1}^{n_{\rm{n}}}{\lambda_{k}}}$, where $\lambda_{i}$ is the $i$'th largest eigenvalue of the data covariance matrix. The total energy plot indicates at least how many reduced particles are required to capture the fullspace information. As seen in Fig. \ref{fig:PCA energy}, 8 reduced particles can cover approximately 95\% of the energy. Fig. \ref{fig:reconst_ex_narrow} and \ref{fig:reconst_ex_wide} show two examples of reconstructed particles from subspace. Although some information is lost, 8 particles can still describe the entire motion of the full system. Thus, we use 8 reduced particles (i.e. PCA modes) in our application. We compute and store the loading matrix $W$ (and mean $\mu$) for reduction and reconstruction.

\subsection{GNN Training}
The GNN model can now be trained using the reduced data. The inputs and outputs are the reduced particle positions and their accelerations, respectively. When computing the particle accelerations, we use finite difference as
\begin{align}
    \dot{p}_{k}\Delta t = {p}_{k} - p_{k-1}, \quad
    \ddot{p}_{k}\Delta t = \dot{p}_{k} - \dot{p}_{k-1},
\end{align}
leading to
\begin{align}
    \ddot{p_k}(\Delta t)^{2} = p_{k} - 2p_{k-1}+p_{k-2}.
\end{align}
We use $\ddot{p_k}(\Delta t)^{2}$ as labels and thus the output of GNN.
Then, the velocity and position update laws with semi-implicit implicit Euler scheme as in \cite{sanchez2020learning};
\begin{align}
    \dot{p}_{k+1} \Delta t = \dot{p}_{k} \Delta t+ \ddot{p}_{k}(\Delta t)^{2}, \quad p_{k+1} = p_{k} + \dot{p}_{k+1}\Delta t.
\end{align}
This way, time step $\Delta t$ is absorbed in the dynamics; thus we do not need to explicitly consider it in the GNN model.

\subsection{DDP States}
In order to perform DDP, the initial and target states are required. For the initial state, we use a rectangular-shaped mass of normal granular particles, which is generated by MPM-NGF and projected into subspace by PCA. The original shape and the reconstructed one are shown in Fig. \ref{fig:schematic} top left. Notice that due to the reconstruction error, the corners are not square. The state of our system needs the first $C+1$ time steps of positions to feed into the GNN model. In order to obtain them as the initial states, we let this shape collapse in the simulation with no control.

The position temporal sequence is concatenated, forming $X_{k:k+C}$ in Eq. \eqref{eq:DDP X state}. The initial box position $b_{k:k+C}$ is set to $[0.5, 0]$ for all the time steps. The target states are shown in Fig. \ref{fig:DDP shapes} (A), which is then projected into subspace. Following the same concatenation procedure, we get the target state $X_{\rm{g}}$ in Eq. \eqref{eq:numerical_example_quad_cost}. Finally, we note that we use $C=5$ as in \cite{sanchez2020learning}.

\subsection{Results}
We perform DDP on the trained GNN model in the subspace. We can then recover the fullspace particles by applying PCA inverse (see Eq. \eqref{eq:PCA reconstruction}). In order to validate our system's dynamics and control commands, we compare not only the reconstructed target and reconstructed terminal states, but also the target and terminal states in fullspace obtained by applying the DDP control commands to the physics-based simulator. Therefore, we move the box in the physics-based simulation following the trajectory of the box from DDP. We perform two experiments with different target states. Where the target shapes have a steep slope (example 1) and a mild slope (example 2).

\begin{figure}[!t]
     \centering
     \begin{subfigure}[b]{0.44\columnwidth}
         \centering
         \includegraphics[trim={0.4cm 0.3cm 0.2cm 0cm},clip,width=\columnwidth]{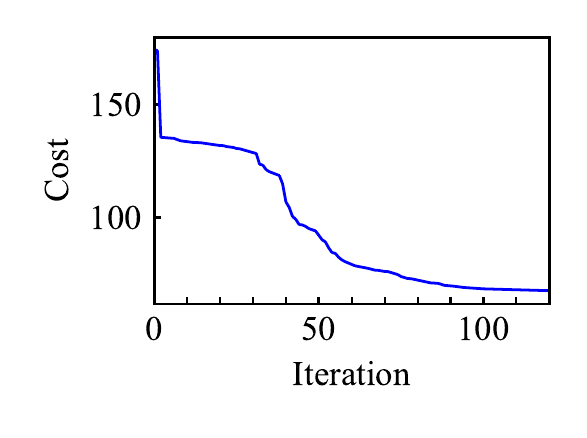}
         \caption{Cost reduction}
         \label{fig:DDP cost reduction}
     \end{subfigure}
     \begin{subfigure}[b]{0.44\columnwidth}
         \centering
         \includegraphics[trim={0.4cm 0.3cm 0.2cm 0cm},clip,width=\columnwidth]{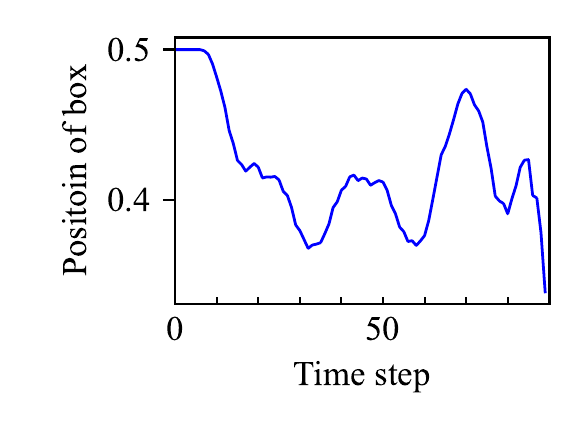}
         \caption{$x$ position of the box}
         \label{fig:DDP box traj}
     \end{subfigure}
     \begin{subfigure}[b]{\linewidth}
         \centering
         \includegraphics[width=\textwidth]{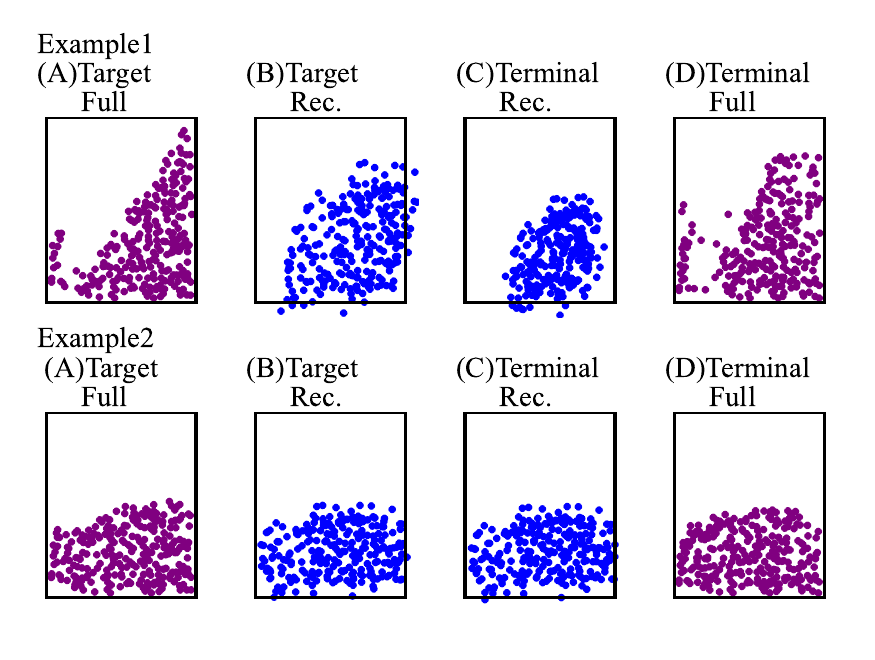}
         \caption{Target and terminal states. Full and reconstructed particles are drawn in purple and blue, respectively. (A)'s are target shapes in fullspace, (B)'s are PCA reconstructed target shapes in fullspace, (C)'s are PCA reconstructed terminal state from GNN model with DDP control commands, and (D)'s are terminal states from physics-based simulator with DDP control commands.}
         \label{fig:DDP shapes}
    \end{subfigure}
    \caption{Model results and validation.}
    \label{fig:DDP result}
\end{figure}

In Fig. \ref{fig:DDP cost reduction}, the cost convergence over DDP iterations in example 1 is presented. It indicates a typical behavior of the cost in DDP (i.e. a quick reduction in the early iterations and a smooth reduction in the end). Also, Fig. \ref{fig:DDP box traj} illustrates the temporal x-positions of the box obtained from DDP.

The visualized results are shown in Fig. \ref{fig:DDP shapes} including two examples with different target shapes. Specifically, (A)'s are the target shapes in fullspace, (B)'s are the PCA reconstructed target shapes in fullspace, (C)'s are the PCA reconstructed terminal state from the GNN model with the DDP control commands, and (D)'s are the terminal states from the physics-based simulator with the DDP control commands. Note the penetrations are due to the reconstruction error. Ideally, the PCA reconstructed (B and C)'s should be exactly the same. In overall, the controller seems to know in which direction it should move the box to achieve the target shape. The root mean square errors (RMSEs) between the fullspace target and terminal states (A and D), and the reconstructed target and terminal states (B and C) are summarized in Table \ref{tab:DDPres} for the two examples.

We deduce that a large portion of the error between (B and C) is due to the DDP method. DDP is a local method that approximates value function around the nominal trajectory; Thus it converges to a local minimum. Tuning the weight matrices might partially resolve this issue. In addition to the dynamics error from GNN, the error between (B and C) also affects the error between (A and D) as they should be equally shaped as well. One resolution we aim to try is to investigate the feedback gain effect. In fact, in the process of applying the computed DDP control commands to the physics-based simulator, the feedback gains $K$ in Eq. \eqref{eq:delta-u-star} can modify the control commands via $K\delta p_{k} = K[p^{\ast}_{k}-W^{\tr}(q_{k}-\mu^{\tr})]$ where $p^{\ast}_{k}$ and $q_{k}$ are the optimal state in the subspace and the current state in the fullspace, respectively. Note that an additional integration is required since the DDP control commands are accelerations. This modification might alleviate the error between (A and D) by keeping (D) closer to (C), which will be tested in the future. For dynamics error, we observe that sometimes DDP commands low amplitude but high frequency control, which is not included in the training dataset. As a result, GNN dynamics shows different behavior as that of physics-based model. This is noticeable in example 2 where error between (A and D) is much larger than that of (B and C). This situation can be avoided by adding state and/or control constraints to the problem using constrained version of DDP, avoiding motion in high frequency \cite{aoyama2021DDP}.

We also speculate that the number of PCA modes (reduced particles) could be another cause of the overall error. Consequently, we tried using more modes with the hope of potentially increasing the accuracy of the entire model. Of course, more modes reduced reconstruction error of PCA. However, it turned out that a higher number of reduced particles makes the optimization problem more challenging with underactuation, i.e., many states with only one dimensional control. With 16 particles (32 states in 2D), DDP gets stuck in a local minima, and thus its cost plateaus without reaching the target shape. It will be further examined in future work.

\begin{table}[!t]
    \footnotesize
    \caption{RMSEs between target and terminal states. Alphabets correspond to those of Fig. \ref{fig:DDP shapes}.}
    \begin{center}
    \begin{tabular}[h]{c|c|c} \toprule
    Example
    & \makecell{Target-terminal MSE\\Fullspace (A and D)\\$\times$ \num{e-3}}
    & \makecell{Target-terminal MSE\\Reconstructed (B and C)\\ $\times$\num{e-3}} \\
    \hline
    1
    & 38.6  
    & 22.1
    \\
    2
    & 17.9
    & 2.20
    \\
    \hline
    \end{tabular}
    \end{center}
    \label{tab:DDPres}
\end{table}

\section{Conclusion}
\label{sec:conclusion}
In this work, to form the granular materials into target shapes, we performed an optimal control using DDP with the PCA-reduced dynamics learned by the GNN model. We used a rectangular rigid box with granular material inside as an example system. After generating training data via a physics-based simulator (MPM-NGF), we reduced the data dimensionality via PCA and trained a GNN model accordingly. We derived the state-space representation of learned reduced dynamics and implemented the DDP method coupled with it. This method computed the temporal control commands that can form granular particles into a target shape. The control sequence was also applied to the physics-based simulator, showing the validity of our approach. 

The future research direction is to tune the learned model and the optimal control algorithm. For the learned model, it is worth trying to utilize non-linear dimensionality reduction methods, such as Auto-encoder and Kernel PCA to improve the performance/accuracy. Note that, here, have only used 8 PCA modes, which seems to be very minimal. For the optimal control algorithm, instead of indexing all particles for the cost computation as we do here, we will utilize distributionally robust cost metrics, such as Chamfer and Earthmover's distances to potentially make the optimization process easier and thus can keep more reduced particles to make reconstruted particles closer to the fullspace particles.



\end{document}